# New Developed Pattern on Number System, New Number Field and Their Applications in Physics

Yi-Fang Chang

Department of Physics, Yunnan University, Kunming, 650091, China

(e-mail: yifangchang1030@hotmail.com)

**Abstract**：Based on a brief review on developments of number system, a new developed pattern is proposed. The quaternion is extended to a matrix form *aI+bC+cB+dA*, in which the unit matrix I and three special matrices C,B,A correspond to number 1 and three units of imaginary number i,j,k, respectively. They form usually a ring. But some fields may be composed of some special 2-rank, even n-rank matrices, for example, three matrices *aI+bC*, *aI+cB*, *aI+dA* and so on. It is a new type of hypercomplex number fields. Finally, the physical applications and possible meaning of the new number system is researched.

**Key words**：number system, matrix, complex number, ring, field, physical meaning

**1.Introduction**

In mathematics first base is number [1]. Basic fields are the field of rational numbers Q, the field of real numbers R and the field of complex numbers C in the number system [2]. Further developments rise and fall. In 1843 W.R.Hamilton obtained a quaternion *a+bi+cj+dk*, here $i^2 = j^2 = k^2 = -1$. But, the quaternion is a ring. He and Maxwell, Heaviside, et al., discussed various applications of the quaternion. Then Graves and Cayley derived biquaternions. Grassmann researched another hypercomplex. In 1873 W.K.Clifford obtained the Clifford number. It is yet other hypercomplex, or the quasi-quaternion $q + \omega Q$, here q and Q are real quaternion, and $\omega^2$ =1. In the Clifford algebra there have elements 1, $e_1, e_2,...e_{n-1}$, here $e_i^2 = -1$, $e_i e_j = -e_j e_i$. In 1889 A.Hurwitz proved that real number, complex number, real quaternion and Clifford quasi-quaternion are only a linear algebra satisfied the multiplication axiom [3]. Then there are the field of algebraic numbers, the field of rational functions, the field of trigonometric functions, the field of meromorphic functions and Hensel p-adic fields and so on [1,2].

In algebra there is Frobenius Theorem: The only division algebras of finite rank over the real number field are itself, the complex number field, and the quaternion algebra [2]. But, in the nonstandard analysis (NSA) there is already new hyper-real number: infinitesimal and infinite [4]. Now the extension of field is an important part of algebra [5]. Here we research that the complex number system is developed to a method by the representation of matrix, and derive a type of the hypercomplex fields. It is a new extended pattern on the number system.

**2.Number system composed by matrix**

Using a matrix representation of quaternion, the complex number system is represented by matrix, here 1→I=$\begin{pmatrix} 1 & 0 \\ 0 & 1 \end{pmatrix}$ is unit matrix, and i→C=$\begin{pmatrix} 0 & -1 \\ 1 & 0 \end{pmatrix}$ and $C^2 = -I$. Assume that a new complex positive number j, here $j^2=1$. Then, let j→B=$\begin{pmatrix} 0 & 1 \\ 1 & 0 \end{pmatrix}$, $B^2=I$, correspond to an inversion operator. k=ij→A=$\begin{pmatrix} -1 & 0 \\ 0 & 1 \end{pmatrix}$=CB, here $A^2 = I$ and ji=-k→-A. A set (1,i,j,k) is another quaternion, and relates the Clifford fourfold quaternion [6,7]. This corresponds to that a field is extended to a ring.

Further, the matrix representation of this number system may extend to higher rank matrix and high dimensional space, and both may be combined each other.

1. We extend to higher unit matrix I, and a main diagonal matrix becomes a secondary diagonal matrix, i.e.,

$$B=\begin{pmatrix} & & 1 \\ & 1 & \\ 1 & & \end{pmatrix} \text{ and so on, to } \begin{pmatrix} & & & 1 \\ & & 1 & \\ & \cdots & & \\ 1 & & & \end{pmatrix}. \quad (1)$$

Moreover, some main diagonal matrices are introduced:

$$A=\begin{pmatrix} -1 & & \\ & -1 & \\ & & 1 \end{pmatrix} \text{ or } \begin{pmatrix} -1 & & \\ & 1 & \\ & & 1 \end{pmatrix} \text{ and so on, to } \begin{pmatrix} -1 & & & \\ & -1 & & \\ & & \cdots & \\ & & & 1 \end{pmatrix}, \quad (2)$$

and these matrices become the secondary diagonal matrices, i.e.,

$$C=\begin{pmatrix} & & -1 \\ & -1 & \\ 1 & & \end{pmatrix} \text{ or } \begin{pmatrix} & & -1 \\ & 1 & \\ 1 & & \end{pmatrix} \text{ and so on, to } \begin{pmatrix} & & & -1 \\ & & -1 & \\ & \cdots & & \\ 1 & & & \end{pmatrix}. \quad (3)$$

Here $B^2=I$, $A^2=I$, only the changes of $\begin{pmatrix} & & -1 \\ & -1 & \\ 1 & & \end{pmatrix}^2 = \begin{pmatrix} -1 & & \\ & 1 & \\ & & -1 \end{pmatrix} = \begin{pmatrix} & & -1 \\ & 1 & \\ 1 & & \end{pmatrix}^2$ and

$\begin{pmatrix} & & & -1 \\ & & -1 & \\ & -1 & & \\ 1 & & & \end{pmatrix}^2 = \begin{pmatrix} -1 & & & \\ & 1 & & \\ & & 1 & \\ & & & -1 \end{pmatrix}$ and so on, are the biggest.

2. Bothe combines each other

$$I+C=\begin{pmatrix} 1 & -1 \\ 1 & 1 \end{pmatrix}, I+B=\begin{pmatrix} 1 & 1 \\ 1 & 1 \end{pmatrix}, I+A=\begin{pmatrix} 0 & 0 \\ 0 & 2 \end{pmatrix}=2I_{22} \text{ and } C+B=\begin{pmatrix} 0 & 0 \\ 2 & 0 \end{pmatrix}=2I_{21}. \quad (4)$$

It is analogue with the elementary matrix $I_{ij}$, in which only these terms on i and j are 1, other elements are 0.

$$(I+C)(I+C)=2C,\ (I+B)(I+B)=2(I+B),\ (I+A)(I+A)=2(I+A). \tag{5}$$

The after two forms are DD=2D $\rightarrow D^n = 2^{n-1}D$. And $(I+C)^3 = 2(C-I)$, $(I+C)^4 = -4I$ and so on. In the quantum field theory [8] the creation operator of particle number is

$$a^+ = \begin{pmatrix} 0 & 0 \\ 1 & 0 \end{pmatrix} = \frac{C+B}{2}, \tag{6}$$

the annihilation operator is

$$a = \begin{pmatrix} 0 & 1 \\ 0 & 0 \end{pmatrix} = \frac{B-C}{2}, \tag{7}$$

and the particle number operator is

$$N = a^+ a = \begin{pmatrix} 0 & 0 \\ 0 & 1 \end{pmatrix} = \frac{I+A}{2},\ aa^+ = \begin{pmatrix} 1 & 0 \\ 0 & 0 \end{pmatrix} = \frac{I-A}{2}. \tag{8}$$

Eqs.(6), (7) and (8) correspond to four numbers. Moreover,

$$a^+ a + aa^+ = \begin{pmatrix} 1 & 0 \\ 0 & 1 \end{pmatrix} = I,\ a^+ + a = \begin{pmatrix} 0 & 1 \\ 1 & 0 \end{pmatrix} = B. \tag{9}$$

Then we may develop to various quaternion and n-element quantity, for example, I+A(I+I+I+…), I+B(I+I+I+…), I+C(I+I+I+…), A+I(I+A+B+C+…), etc.

3. A general main diagonal matrix is

$$A = \begin{pmatrix} a_1 & & & \\ & a_2 & & \\ & & \ldots & \\ & & & a_n \end{pmatrix}, \tag{10}$$

and a secondary diagonal matrix is

$$B = \begin{pmatrix} & & & b_1 \\ & & b_2 & \\ & \ldots & & \\ b_n & & & \end{pmatrix}. \tag{11}$$

Their products are still main diagonal matrices or secondary diagonal matrices:

$$AA' = \begin{pmatrix} a_1 a'_1 & & & \\ & a_2 a'_2 & & \\ & & \ldots & \\ & & & a_n a'_n \end{pmatrix},\ AB = \begin{pmatrix} & & & a_1 b_1 \\ & & a_2 b_2 & \\ & \ldots & & \\ a_n b_n & & & \end{pmatrix},$$

$$BA = \begin{pmatrix} & & & b_1 a_n \\ & & b_2 a_{n-1} & \\ & \ldots & & \\ b_n a_1 & & & \end{pmatrix}, \quad BB' = \begin{pmatrix} b_1 b'_n & & & \\ & b_2 b'_{n-1} & & \\ & & \ldots & \\ & & & b_n b'_1 \end{pmatrix},$$

$$B'B = \begin{pmatrix} b_n b'_1 & & & \\ & b_{n-1} b'_2 & & \\ & & \ldots & \\ & & & b_1 b'_n \end{pmatrix}, \text{ (here } a_i, b_i \in C\text{).} \tag{12}$$

**3. Field composed by matrices**

The matrix $\{aI + bC + cB + dA \mid a,b,c,d \in C\}$ forms a ring for addition and multiplication. Here

$$CB=A, BC=-A; AC=B, CA=-B; AB=C, BA=-C. \tag{13}$$

The multiplication of matrix is non-commutation, and is anticommutation. Since A, B and C have the inverse-elements, the division for matrix is defined by the inverse matrix. Such it is a division ring of real quaternions. But, three matrices $aI+bC$, $aI+cB$ and $aI+dA$ all have inverse-elements, and obey the commutative law and the associative law of multiplication, and construct three similar commutative quotient ring——field. An example is:

$$aI + bC = \begin{pmatrix} a + a_1 i & -b - b_1 i \\ b + b_1 i & a + a_1 i \end{pmatrix}. \tag{14}$$

In usual case the matrix does not obey the commutative law of multiplication, but the main diagonal matrices obey the commutative law. For the n-rank matrices (12), when $a_1 = a_n$, $a_2 = a_{n-1}$,…i.e., AB=BA, and $b_1 b'_n = b_n b'_1$, $b_2 b'_{n-1} = b_{n-1} b'_2$, …i.e., $BB' = B'B$, it obeys the commutative law. In these cases a new number field may be constructed.

In more general case, for example, some special 2-rank matrices are

$$J = \begin{pmatrix} bh + d & b \\ kb & d \end{pmatrix}, \text{ (here h and k are constant).} \tag{15}$$

It obeys the commutative law of multiplication. Therefore, $\{aI + cJ \mid a,c \in C\}$ constructs the field. These matrices include: 1. A supertriangle matrix (for k=0), which may become a main diagonal matrix. 2. $J = b\begin{pmatrix} h & 1 \\ k & 0 \end{pmatrix}$ (for d=0). 3. $J = \begin{pmatrix} d & b \\ kb & d \end{pmatrix}$ (for h=0). Further, this method may be extended to 3-rank, n-rank special matrices. These matrices have unit matrices I, zero element, negative element and inverse element, and addition, etc., and they obey yet commutative law of multiplication, so they may construct a number field. Any rank diagonal matrices all obey the commutative law of multiplication, and may construct the field.

Usual matrix obeys the algorithm of ring. This pattern may construct a big type of general ring by matrix representation, for instance, algebra of Pauli ring and Dirac ring [6]. Moreover, they may construct the field by some special matrices under the certain conditions. For example, the 3-rank matrices:

$$I = \begin{pmatrix} 1 & & \\ & 1 & \\ & & 1 \end{pmatrix}, \quad D_1 = \begin{pmatrix} & 1 & \\ & & 1 \\ 1 & & \end{pmatrix}, \quad D_2 = \begin{pmatrix} 1 & & \\ & -1 & \\ & & 1 \end{pmatrix}, \quad D_3 = \begin{pmatrix} & & 1 \\ & -1 & \\ 1 & & \end{pmatrix}. \quad (16)$$

Because $D_1 D_2 = D_2 D_1 = D_3$, $D_1 D_3 = D_3 D_1 = D_2$ and $D_2 D_3 = D_3 D_2 = D_1$, the matrices of quaternion form $aI + bD_1 + cD_2 + dD_3$, and in which three dual-parts $aI + bD_1$, $aI + cD_2$ and $aI + dD_3$ all construct the general fields. Simultaneously,

$$I = \begin{pmatrix} 1 & & \\ & 1 & \\ & & 1 \end{pmatrix}, \quad D_2 = \begin{pmatrix} 1 & & \\ & -1 & \\ & & 1 \end{pmatrix}, \quad D_4 = \begin{pmatrix} & & 1 \\ & 1 & \\ -1 & & \end{pmatrix}, \quad D_5 = \begin{pmatrix} & & 1 \\ & -1 & \\ -1 & & \end{pmatrix}. \quad (17)$$

Because the products of $D_2, D_4, D_5$ all obey the commutative law of multiplication, a new hypercomplex number

$$aI + bD_2 + cD_4 + dD_5, \quad (18)$$

and in which some dual-parts construct yet the fields. According to the same composing rules, a matrix of quaternion form is composed of these 4-rank matrices:

$$I = \begin{pmatrix} 1 & & & \\ & 1 & & \\ & & 1 & \\ & & & 1 \end{pmatrix}, G_1 = \begin{pmatrix} -1 & & & \\ & 1 & & \\ & & 1 & \\ & & & -1 \end{pmatrix}, G_2 = \begin{pmatrix} & & & 1 \\ & & 1 & \\ & 1 & & \\ 1 & & & \end{pmatrix}, G_3 = \begin{pmatrix} & & & -1 \\ & & 1 & \\ & 1 & & \\ -1 & & & \end{pmatrix},$$

which construct still a field. It is namely that a field is constructed by the matrix as follows:

$$aI + bG_1 + cG_2 + dG_3 = \begin{bmatrix} a-b & & & c-d \\ & a+b & c+d & \\ & c+d & a+b & \\ c-d & & & a-b \end{bmatrix}. \quad (19)$$

For more general cases the eight matrices of 5-rank I and $G_i$ (i=1,2,3,4,5,6,7) are composed of the main diagonal elements, which are respectively (1,1,1,1,1), (-1,1,1,1,-1), (-1,-1,1.-1,-1), (1,-1,1,-1,1) and corresponding secondary diagonal elements. This form

$$a_1 I + a_2 G_1 + a_3 G_2 + a_4 G_3 + a_5 G_4 + a_6 G_5 + a_7 G_6 + a_8 G_7, \text{ (here } a_i \in C\text{)}, \quad (20)$$

constructs yet a field. The eight matrices of 6-rank are composed of the main diagonal elements, which are respectively (1,1,1,1,1,1), (-1,1,1,1,1,-1), (-1,-1,1,1,-1,-1), (1,-1,1,1,-1,1) and corresponding secondary diagonal elements. It constructs still a field.

The twelve matrices of 7-rank are composed of the main diagonal elements, which are respectively (1,1,1,1,1,1,1), (-1,1,1,1,1,1,-1), (-1,-1,1,1,1,-1,-1), (-1,-1,-1,1,-1,-1,-1), (1,-1,-1,1,-1,-1,1), (1,-1,1,-1,1,-1,1) and corresponding secondary diagonal elements. It constructs a field. The twelve matrices of 8-rank are composed of the main diagonal elements, which are respectively (1,1,1,1,1,1,1,1), (-1,1,1,1,1,1,1,-1), (-1,-1,1,1,1,1,-1,-1), (-1,-1,-1,1,1,-1,-1,-1), (1,-1,-1,1,1,-1,-1,1), (1,-1,1,-1,-1,1,-1,1) and corresponding secondary diagonal elements. It constructs yet a field and so on. These fields are composed of even elements. In a word, the number system extends to matrix, which may construct very abundant various fields.

**4. Physical applications and possible meaning of new number system**

The multiplication of quaternion is non-commutative, which may describe the composition of rotation on the rigid body. In relativity if the four dimensional space-time corresponds to a quaternion, it should be that three imaginary numbers correspond to three dimensional space, and a real number corresponds to one dimensional time. Assume that a four dimensional space-time is:

$$(Ix)^2 + (By)^2 + (Cz)^2 + (Act)^2. \tag{21}$$

This is a new representation of four dimensional space-time. A similar higher rank matrix corresponds to a metric matrix $\begin{pmatrix} 1 & & & \\ & 1 & & \\ & & 1 & \\ & & & -1 \end{pmatrix}$ of four dimensional space-time of special relativity. It becomes a light-cone for two dimensional space-time, and four dimensional geometry corresponds to a hyperboloid of two sheets. The metric matrix of general relativity corresponds to a general 4-rank matrix, and hyper-real number. Segel discussed the infinite and the infinitesimal in models for natural phenomena [9].

In quantum mechanics matrix is a basic tool, and q-number is non-commutative. They correspond to operator and ring. Therefore, the hypercomplex system (1,i,j,k) relates three dimensional Pauli matrices in quantum mechanics [10], i.e.,

$$\sigma_x = B, \quad \sigma_y = i\begin{pmatrix} 0 & -1 \\ 1 & 0 \end{pmatrix} = iC, \quad \sigma_z = -A. \tag{22}$$

They and I construct four dimensional space-time. Further, this relates the Dirac matrices and extensions, in which $G_3$ is namely $\gamma_2$ in Dirac matrices. Any products of Pauli and Dirac matrices all are closed. 1 and Pauli matrices are the basic elements of a ring, which forms the Pauli algebra. Dirac matrices construct another ring, and correspond to another algebra and Clifford number [6]. 1i=i1 is commutative, and corresponds to boson; ij=-ji is anticommutative, and corresponds to fermion. The complex number system relates the supersymmetry and Graded Lie Algebras (GLA), in which the commutation rules for the (D+d) dimensional GLA [11] are:

$$[A_m, A_n] = A_m A_n - j^2 A_n A_m = f_{mn}^l A_l. \tag{23}$$

Even elements and Lie Algebra corresponds to boson and integral spin.

$$\{Q_\alpha, Q_\beta\} = Q_\alpha Q_\beta - i^2 Q_\beta Q_\alpha = F_{\alpha\beta}^m A_m. \tag{24}$$

Odd elements and GLA corresponds to fermion and half-integral spin. The unification of A and Q corresponds to unification of i and j. Moreover,

$$[A_m, Q_\alpha] = A_m Q_\alpha - k^2 Q_\alpha A_m = S_{m\alpha}^\beta Q_\beta. \tag{25}$$

The above development of complex number system will redound to various unified representations between Lie Algebra and GLA, between commutation and anticommutation relations, between boson and fermion, and between Bose-Einstein statistics and Fermi-Dirac statistics. This relates the supersymmetry, and a unified statistics and possible violation of Pauli exclusion principle [12-15]. Recently D.G.Pavlov discussed the relations among the hypercomplex numbers, the associated metric space and extension of the relativistic hyperboloid [16].

In special relativity the Lorentz transformation may be represented as a matrix form. The supersymmetric change between bosons and fermions may apply matrix. Combining both, we proposed a mathematical physical law: Bosons correspond to real number, and fermions correspond to imaginary number. Such bosons and fermions consist of even and odd fermions, respectively, which just corresponds to even and odd imaginary numbers are real and imaginary number [17,18]. Combining the development of complex number, we developed the complex number from the fractal dimension, and researched its applications in mathematics and physics, and the fractal space-time theory [19,20]. From this it may be constructed that the higher dimensional, fractal, complex and hypercomplex space-time theory covers all [19,20].

**5.Conclusion**

1.The number extends to a matrix representation, i.e., number 1 and imaginary number units (i,j,k) correspond to a unit-matrix I and the special matrices (C,B,A). 2.The matrix *aI+bC+cB+dA* constructs a ring, and corresponds to the forms of quaternion. 3.Three matrices *aI+bC*, *aI+cB* and *aI+dA* all are fields. 4.We may extend some special 2-rank, even higher-rank matrices, and construct some fields. 5.This new developed pattern of number system is possibly some physical meaning. We predict that if the field can be extended, this will be able to apply to many regions.